\documentclass[12pt]{article}
\usepackage{a4,latexsym,amsmath,amssymb,amsthm,enumerate,eucal,rawfonts,
  graphicx,xcolor,subfig}
\newtheorem{theorem}{Theorem}
\newtheorem{lemma}[theorem]{Lemma}
\newtheorem{proposition}[theorem]{Proposition}
\newtheorem{claim}{Claim}

\newtheorem{problem}[theorem]{Problem}

\newcommand{\claimproofend}{\hspace*{.1mm}\hspace{\fill}$\diamond$}
\newenvironment{claimproof}{}{\claimproofend\par\vspace{2mm}}

\newcommand\Setx[1] {\left\{{#1}\right\}}

\newcommand\size[1] {\left|{#1}\right|}

\newcommand{\rconn}[1]{\mathit{rc}(#1)}
\newcommand{\im}[1]{\mathit{im}\ #1}
\newcommand{\eps}{\varepsilon}

\newcommand{\Hstar}{H^*}
\newcommand{\cstar}{c^*}
\newcommand{\sm}{\setminus}

\newtheorem{xcasehdr}{Case}
\newtheorem{xxcasehdr}{Subcase}[xcasehdr]

\newenvironment{xcase}[1]%
{\vspace{-1mm}\par\noindent
  \xcasehdr{#1}\upshape
  \vspace{2mm}\par\noindent}
{\hspace*{0mm}\hspace{\fill}$\diamond$\par\vspace{2mm}}

\newlength{\zero}
{
  \par\noindent
  \hangindent\parindent
  \hangafter=0
  \xxcasehdr{#1}\vspace{2mm}\par\upshape\noindent}
{\hspace*{0mm}\hspace{\fill}$\triangle$\par\vspace{4mm}}

\newcommand{\fig}[1]{\includegraphics[page=#1]{rainbow-figs.pdf}}%
\newcommand{\sfig}[2]{\subfloat[#2]{\fig{#1}}}%

\newcommand{\sfigtop}[2]{\newbox{\pic}\sbox{\pic}{\fig{#1}}%
  \subfloat[#2]{\vbox to\ht\base{\hbox to\wd\pic{\usebox\pic}}}}
\newcommand{\hf}{\hspace*{0pt}\hspace{\fill}\hspace*{0pt}}

\date{\today}
\title{\textbf{The rainbow connection number\\of 2-connected graphs}}

\author{Jan Ekstein$^{\,\ast}$ \and %
  P\v{r}emysl Holub$^{\,\ast}$ \and %
  Tom\' a\v s Kaiser$^{\,\ast}$ \and %
  Maria Koch$^{\,\dagger}$ \and %
  Stephan Matos Camacho$^{\,\dagger}$ \and %
  Zden\v ek Ryj\'{a}\v{c}ek\footnote{Department of Mathematics and
    Institute for Theoretical Computer Science, University of West
    Bohemia, Pilsen, Univerzitn\'{\i} 8, 306 14 Pilsen, Czech
    Republic. E-mail: \texttt{\{ekstein, holubpre, kaisert,
      ryjacek\}@kma.zcu.cz}. Research supported by grants MEB 101014
    and 1M0545 of the Czech Ministry of Education.} \and %
  Ingo Schiermeyer\footnote{Institut f\"{u}r Diskrete Mathematik und
    Algebra, Technische Universit\"{a}t Bergakademie Freiberg, 09596
    Freiberg, Germany. E-mail: \texttt{\{Maria.Koch, Matos,
      Ingo.Schiermeyer\}@math.tu-freiberg.de}. Research supported by
    the DAAD-PPP project ``Colourings, cycles and closures''.}}

\date{}

\begin{document}

\maketitle

\begin{abstract}
  The rainbow connection number of a graph $G$ is the least number of
  colours in a (not necessarily proper) edge-colouring of $G$ such
  that every two vertices are joined by a path which contains no
  colour twice. Improving a result of Caro et al., we prove that the
  rainbow connection number of every 2-connected graph with $n$
  vertices is at most $\lceil n/2 \rceil$. The bound is optimal.
\end{abstract}

\section{Introduction}
\label{sec:intro}

We investigate a problem related to the concept of rainbow connection
in graphs, introduced by Chartrand et al.~\cite{CJMZ}. Let $G$ be an
undirected graph with a colouring $c$ of the edges, which is not
assumed to be proper (that is, adjacent edges may get the same
colour). A subgraph $H$ of $G$ is \emph{rainbow} (with respect to $c$)
if no two edges of $H$ have the same colour under $c$. The
edge-coloured graph $(G,c)$ is \emph{rainbow-connected} if every pair
of vertices is joined by a rainbow path. The \emph{rainbow connection
  number} of $G$, denoted by $\rconn G$, is the least number of
colours in a colouring which makes $G$ rainbow-connected.

If $G$ has $n$ vertices, then $\rconn G\leq n-1$ as one may colour
each edge of a spanning tree of $G$ with a different colour, and use
one of these colours for all the remaining edges. Chartrand et
al.~\cite{CJMZ} determined the rainbow connection number of several
classes of graphs, such as the complete multipartite graphs. The
rainbow connection number has been studied for further graph classes
in \cite{CLRTY} and for graphs with fixed minimum degree in
\cite{CLRTY, KY, S1, S2}.

The computational complexity of rainbow connectivity has been studied
in \cite{CFMY} where it is proved that determining the rainbow
connection number is an NP-complete problem. Indeed, it is already
NP-complete to decide whether $\rconn G$ equals two~\cite{CFMY}. More
generally, it was shown in \cite{LT} that for any fixed $k \geq 2$,
deciding if $\rconn G = k$ is NP-complete.

Caro et al.~\cite{CLRTY} proved the following upper bound for the
rainbow connection number of a 2-connected graph:
\begin{theorem}[\cite{CLRTY}]
  If $G$ is a 2-connected graph on $n$ vertices, then 
  \begin{equation*}
    \rconn G \leq \frac n2 + O(\sqrt{n}).
  \end{equation*}
\end{theorem}

In this paper, we improve this upper bound to an optimal one, which is
attained, e.g., for all odd cycles:
\begin{theorem}
  \label{t:main}
  For any graph $G$ with $n$ vertices,
  \begin{equation*}
    \rconn G \leq \Bigl\lceil\frac n2 \Bigr\rceil.
  \end{equation*}
\end{theorem}

The proof of Theorem~\ref{t:main} will be presented in
Section~\ref{sec:proof-theorem}. It is based on several lemmas which
are established in the Sections~\ref{sec:paths} and
\ref{sec:extending}.


In the remainder of this section, we fix the necessary notation and
terminology. For the terms not defined here, as well as for broader
background, the reader may wish to consult~\cite{BM}.

Our graphs are finite, undirected and simple; in particular, parallel
edges are not allowed. The vertex and edge sets of $G$ are denoted by
$V(G)$ and $E(G)$. The number of vertices of a graph $G$ is denoted by
$\size G$. A path with endvertices $x$ and $y$ is referred to as an
\emph{$xy$-path}. For $H\subseteq G$, an \emph{$H$-path} is a path
disjoint from $H$ except for its endvertices, which are contained in
$H$. If $P$ is a path in $G$ and $u,v\in V(P)$, then $uPv$ denotes the
unique subpath of $P$ with endvertices $u$ and $v$. If $Q$ is a path
with $v,w\in V(Q)$, then $uPvQw$ denotes the concatenation of $uPv$
and $vQw$. This may in general be a walk rather than a path, but this
distinction is not too important since any rainbow $uw$-walk contains
a rainbow $uw$-path.

Throughout this paper, the term \emph{colouring} will be used as in
this section --- meaning an edge-colouring which is not necessarily
proper. It is convenient to call a colouring $c$ of $G$
\emph{rainbow-connecting} if $(G,c)$ is rainbow-connected. Since a
colouring $c$ is, formally, a function defined on $E(G)$, it makes
sense to let $\im c$ denote the set of all colours used by $c$.

A somewhat technical strengthening of the concept of a
rainbow-connecting colouring will be useful in our proofs. Let us call
a rainbow path in $(G,c)$ \emph{blocking} if it uses all colours in
$\im c$. Given vertices $x,y\in V(G)$, we say that $y$ is
\emph{blocked for $x$} if all rainbow $xy$-paths are blocking. A
colouring is \emph{safe} if for each vertex $x$, there is at most one
blocked vertex $y$. A colouring which is both safe and
rainbow-connecting is said to be \emph{safely rainbow-connecting}.

Let $A$ be a set of colours and let $c_0$ be a colouring of $G$. A
subgraph $H\subseteq G$ is \emph{$A$-free} in $(G,c_0)$ if no colour
from $A$ is used by $c_0$ on an edge of $H$. To simplify the notation,
we abbreviate, e.g., `$\Setx{\alpha,\gamma}$-free' to
`$\alpha\gamma$-free'.

In the proofs in this paper, we use the symbol $\diamond$ to mark the
end of the proof of a claim. The same symbol is used at the end of the
discussion of each case in a case analysis.


\section{A lemma on paths}
\label{sec:paths}

In this section, we prove a lemma which is one of the key parts of our
argument.

Let $H$ be a subgraph of $G$. A subgraph $H'$, $H\subseteq H'\subseteq
G$, is a \emph{$k$-extension} of $H$ with \emph{path sequence}
$(P_1,\dots,P_k)$ if each $P_i$ is a $(H\cup P_1\cup \dots
P_{i-1})$-path with at least one endvertex in $V(P_1\cup\dots\cup
P_{i-1})$, and
\begin{equation*}
  H' = H\cup P_1\cup\dots\cup P_k.
\end{equation*}
If $k$ is not important, we just say that $H'$ is an extension of
$H$. An extension is \emph{even} (\emph{odd}) if all the paths in the
path sequence are even (odd, respectively). (Recall that even paths
are those with an even number of edges.)

In the proof of the lemma, we will need a concept similar to that of
an $H$-bridge as introduced by Tutte (see,
e.g.,~\cite[Section~9.4]{BM}). A \emph{weak $H$-bridge} is any
component $B$ of $G-E(H)$ containing at least one edge. The vertices
of $B\cap H$ are called the \emph{attachment vertices} of $B$.

\begin{lemma}\label{l:paths}%
  Let $H$ be a connected subgraph of a 2-connected graph $G$ such that
  the following holds:
  \begin{enumerate}[\quad(i)]
  \item all $1$-extensions of $H$ are even,
  \item there is no even $2$-extension of $H$.
  \end{enumerate}
  Then for any $k\geq 1$ and any extension of $H$ with path sequence
  $(P_1,\dots,P_k)$, $P_1$ is even and all the other paths $P_i$ are
  odd.
\end{lemma}
\begin{proof}
  Suppose $H' = H\cup P_1\cup\dots\cup P_k$ is a $k$-extension of $H$,
  and proceed by induction on $k$. The case $k \leq 2$ follows from
  the assumptions, so we assume that $k > 2$, and that the assertion
  holds for $j$-extensions of $H$ with $j < k$. In particular, each
  $P_i$ ($2\leq i < k$) is odd, and therefore is not an $H$-path by
  condition (i). Observe also that we may assume $\size H \geq 2$,
  since otherwise there exists no $H$-path and the statement is
  trivially true.
  
  For the sake of a contradiction, suppose that $P_k$ is an even
  $(H\cup P_1\cup \dots \cup P_{k-1})$-path with at least one
  endvertex in $P_1\cup \dots \cup P_{k-1}$. We will find either an
  even 2-extension of $H$ (contradicting (ii)), or an extension of $H$
  with a shorter path sequence terminating with $P_k$ (contradicting
  the induction hypothesis).

  Let $G/H$ be the multigraph obtained by contracting all the edges of
  $H$. Since $H$ is connected, the contraction merges all the vertices
  of $H$ into one vertex $*_H$. 

  \begin{claim}\label{cl:paths-bi}
    The graph $G/H$ is bipartite.
  \end{claim}
  \begin{claimproof}
    Suppose, for the sake of a contradiction, that $C_*$ is an odd
    cycle in $G/H$. If $*_H\notin V(C_*)$, then $C_*$ is also a
    subgraph of $G$, and we can find two vertex-disjoint paths, each
    joining a vertex of $H$ to a vertex of $C_*$ and having no
    internal vertices in $H\cup C_*$. Combining the paths with a
    suitable subpath of $C_*$, we obtain an odd $H$-path, in violation
    of (i).

    Thus, $C_*$ must contain $*_H$. The subgraph of $G$ corresponding
    to $C_*$ is either a path or a cycle. If it is a path, then it is
    an odd $H$-path violating condition (i). Hence, $G$ contains an
    odd cycle $C$ containing exactly one vertex $u_1$ of $H$. Using
    the 2-connectedness of $G$ and the assumption that $\size H\geq
    2$, we can find a path which joins a vertex $u_2\neq u_1$ of $H$
    to a vertex of $C$, and has no internal vertices in $H\cup C$. The
    concatenation of this path with a suitable subpath of $C$ is an
    odd $H$-path, a contradiction.
  \end{claimproof}

  In view of Claim~\ref{cl:paths-bi}, we can let $b_*$ be a
  2-colouring of the vertices of $G/H$ which is proper (adjacent
  vertices get different colours). Consider the corresponding
  2-colouring $b$ of $G-E(H)$ obtained by assigning each vertex $w$
  the colour $b_*(w)$ if $w\notin V(H)$ and $b_*(*_H)$ otherwise.

  Let $B$ be the weak $H\cup P_1$-bridge of $H'$ containing $P_k$. We
  let $A$ denote the set of attachment vertices of $B$ which are
  contained in $P_1$ (that is, $A=B\cap P_1)$. Furthermore, we let
  $J\subseteq\Setx{2,\dots,k}$ be the set of indices $i$ such that
  $P_i$ is contained in $B$. We write $J = \Setx{i_1,\dots,i_\ell}$
  with $i_1 < \dots < i_\ell$ and note that $i_\ell = k$.

  Clearly, if $\ell = 1$ (that is, if $B = P_k$), then $P_k$ is odd as
  we would otherwise obtain an even 2-extension of $H$. In the sequel,
  we will therefore assume that $\ell \geq 2$.

  \begin{claim}\label{cl:paths-diff}
    Any two vertices of $A$ have different colours under $b$.
  \end{claim}
  \begin{claimproof}
    Let $x,y\in A$. Suppose that $b(x) = b(y)$. Since $B$ is
    connected, it contains an $xy$-path $P$. As $P$ is edge-disjoint
    from $H$, the colours of the vertices along $P$ alternate and
    hence $P$ is even. Adding $P$ to $P_1$, we obtain an even
    2-extension of $H$, a contradiction with (ii).
  \end{claimproof}

  By Claim~\ref{cl:paths-diff}, $A$ contains at most two vertices.

  \begin{claim}\label{cl:paths-one}
    The size of $A$ is at most one.
  \end{claim}
  \begin{claimproof}
    For the sake of a contradiction, suppose that $A = \Setx{x,y}$. We
    have $b(x)\neq b(y)$, so if $B\cap H$ contains any other vertex
    $z$, then the colour of $z$ matches that of $x$ or $y$, and
    consequently $B$ contains either an even $xz$-path or an even
    $yz$-path. This provides us again with an even 2-extension of $H$,
    contradicting (ii).

    Thus, $x$ and $y$ are the only two attachment vertices of
    $B$. Necessarily, $P_{i_1}$ is an $xy$-path. Let $Q$ be the unique
    $H$-path obtained by concatenating $P_{i_1}$ with subpaths of
    $P_1$. Note that for any $j$ ($2\leq j \leq \ell$) and any
    endvertex $z$ of $P_{i_j}$ with $z\in V(P_1)$, we also have $z\in
    V(Q)$. Therefore, $H\cup Q\cup B$ can be obtained as the
    $\ell$-extension of $H$ with path sequence
    \begin{equation*}
      (Q, P_{i_2},\dots,P_{i_\ell}).
    \end{equation*}
    Since $2\leq \ell \leq k-1$, the induction hypothesis implies that
    the path $P_{i_\ell}$ is odd. Furthermore (again since we assume
    that $\ell \geq 2$), $P_{i_\ell}$ coincides with $P_k$. This is a
    contradiction with the assumption that $P_k$ is even.
  \end{claimproof}

  We cannot have $A=\emptyset$, since at least one endvertex of
  $P_{i_1}$ is required to lie on $P_1\cup\dots\cup P_{i_1-1}$ and
  thus (by the choice of $i_1$) it must actually be contained in
  $P_1$. 

  Hence, $A$ contains a single vertex, say $A=\Setx x$. The argument
  for this case is similar to that in the proof of
  Claim~\ref{cl:paths-one}. We note that one endvertex of $P_{i_1}$ is
  $x$ and the other endvertex is in $H$. Let $Q$ be an $H$-path
  obtained by concatenating $P_{i_1}$ with a subpath $R$ of $P_1$; of
  the two possibilities for $R$, we choose one where the endvertex of
  $R$ in $H$ is different from $x$. As before, if $z\in V(P_{i_j} \cap
  P_1)$ ($2\leq j \leq\ell$), then $z\in V(Q)$. Consequently, $H\cup
  Q\cup B$ is an $\ell$-extension of $H$ with path sequence
  \begin{equation*}
    (Q, P_{i_2},\dots,P_{i_\ell}).
  \end{equation*}
  Again, we infer from the induction hypothesis that $P_{i_\ell}$ is
  odd, which contradicts the assumption that $P_k = P_{i_\ell}$ is
  even.
\end{proof}


\section{Extending the colourings}
\label{sec:extending}

Throughout this section, let $H$ be a connected subgraph of a graph
$G$, $\size H\geq 3$, and let $c$ be a rainbow-connecting colouring of
$H$. We introduce `standard' ways to extend $c$ to an odd 1-extension
of $H$ and to an even 2-extension of $H$. In this and the following
section, the symbol $\gamma$ will denote a fixed colour which is
assumed to be contained in $\im c$. We make the assumption that
$\size{\im c} \geq 2$.

Suppose first that $P$ is an odd $H$-path, say of length $2k+1$. A
\emph{continuation} of $c$ to $H\cup P$ is any colouring $c'$ which
agrees with $c$ on $H$ and assigns to the edges of $P$, in some
direction, the colours
\begin{equation*}
  a_1,a_2,\dots,a_k,\gamma,a_1,a_2,\dots,a_k,
\end{equation*}
where the $a_i$ ($1\leq i \leq k$) are some distinct colours not
contained in $\im c$. Thus, a continuation of $c$ to $H\cup P$ is not
uniquely determined, but any two such continuations are isomorphic in
the obvious sense, so we may regard them as identical. Note that $c'$
uses $k$ colours not contained in $\im c$, which is half the number of
vertices in $V(P)\sm V(H)$.

\begin{lemma}\label{l:odd}
  Let $c$ be a safely rainbow-connecting colouring of $H$. If $P$ is
  an $H$-path of odd length, then any continuation of $c$ to $H\cup P$
  is safely rainbow-connecting.
\end{lemma}
\begin{proof}
  Let $c'$ be a continuation of $c$ to $H\cup P$. Suppose that $P$ is
  an $H$-path of length $2k+1$ with endvertices $u$ and $v$. We may
  assume that $k\geq 1$ since the statement is trivially true for $k =
  0$.

  To see that $c'$ is rainbow-connecting, we need to exhibit a rainbow
  $xy$-path $R_{xy}$ for each pair $x,y\in V(H\cup P)$. If $x,y\in
  V(H)$, we define $R_{xy}$ as a rainbow $xy$-path in $(H,c)$ which
  exists by the assumption, and we choose a path which is non-blocking
  under $c$ if possible. Note that $R_{xy}$ is rainbow and
  non-blocking under $c'$.

  From now on, we assume that $x\in V(P)$. Suppose first that $y\notin
  V(P)$. Since $P$ is odd, the subpaths $xPu$, $xPv$ are not of the
  same length; without loss of generality, let $xPu$ be the shorter
  one. Note that on $xPu$, $c'$ uses no colour from $\im c$. Thus, if
  we let $R_{xy} = xPuR_{uy}y$ (where $R_{uy}$ has been defined above
  because $u\in V(H)$), then $R_{xy}$ is rainbow in $H\cup P$.

  Assume next that $y\in V(P)$. We may assume that $y\in V(xPv)$. If
  $xPy$ is not rainbow, then it includes a pair of edges with the same
  colour, in which case each colour used by $c'$ on $P$ (including
  $\gamma$) must appear on $xPy$. But then $c'$ is rainbow on $uPx\cup
  vPy$ and uses no colour from $\im c$. Thus, $R_{xy} := xPuR_{uv}vPy$
  is rainbow.

  It remains to show that $c'$ is safe. Let $x$ be a vertex of $H\cup
  P$; we show that at most one vertex is blocked for $x$ under
  $c'$. Since $k > 1$, we may choose a colour $\eps\neq\gamma$ used by
  $c'$ on $P$.

  Suppose that $x\in V(H)$. Since $\eps$ is not used by $c'$ on $H$,
  no $y\in V(H)$ is blocked for $x$ under $c'$.  As for $y\in V(P)$,
  $y$ will only be blocked if $R_{xy}\cap P$ includes all the colours
  used by $c'$ on $P$, possibly except for $\gamma$. This happens only
  if $y$ is incident with the central edge $e$ of $P$. Let $u'$ and
  $v'$ be the endvertices of $e$, where $u'$ is closer to $u$ than to
  $v$. Note that the path $R_{xu'}$ contains $R_{xu}$ as a subpath,
  and similarly $R_{xv} \subseteq R_{xv'}$. Since only one of $u$ and
  $v$ can be blocked for $x$ under $c$, we may assume that $R_{xu}$ is
  not blocking under $c$. Therefore, $R_{xu'}$ is not blocking under
  $c'$ since no colour from $\im c$ is used on the complementary
  subpath $uPu'$ of $R_{xu'}$. Hence, $v'$ is the only vertex which
  may be blocked for $x$ under $c'$.

  The last case to consider is that both $x$ and $y$ are internal
  vertices of $P$. If $R_{xy} \subseteq P$, then the only colour it
  uses from $\im c$ (if any) is $\gamma$. By the assumption that
  $\size{\im c} \geq 2$ (made at the beginning of this section),
  $R_{xy}$ is not blocking. Thus, we may assume that $R_{xy} =
  xPuR_{uv}Py$. There is only one vertex $y$ for which $xPu$ and $vPy$
  cover $\im{c'}\sm \im c$, namely the other vertex of $P$ which is
  incident with edges of the same colours as $x$. Thus, we have shown
  that $c'$ is safe, and the proof is complete.
\end{proof}

Next, we define a continuation of $c$ to an even 2-extension $H\cup
Q\cup Q'$ of $H$, where the length of $Q$ is $2\ell$ and the length of
$Q'$ is $2\ell'$. Suppose that the vertices of $Q$ are
$u_0,\dots,u_{2\ell}$ and the vertices of $Q'$ are
$u'_0,\dots,u'_{2\ell'}$. Let us write $u = u_0$, $v = u_{2\ell}$, $u'
= u'_0$ and $v' = u'_{2\ell'}$. We may assume that the distance of
$u'$ from $H$ in $H\cup Q$ is greater than or equal to the distance of
$v'$ from $H$. In particular, $u'\in V(Q)$. We may also assume that
$u' = u_k$ with $k\leq\ell$, and if $u' = u_\ell$, then $v'\in
V(u_\ell Qv\cup H)$.

We colour the edges of $Q$, in order from $u$ to $v$, as
\begin{equation*}
  a_1,a_2,\dots,a_{\ell-1},a_\ell,\gamma,a_1,a_2,\dots,a_{\ell-1}.
\end{equation*}
The edges of $Q'$, in order from $u'$ to $v'$, will be coloured as
\begin{equation*}
  a_{\ell+1},a_{\ell+2},\dots,a_{\ell+\ell'-1},\gamma,a_\ell,
  a_{\ell+1},a_{\ell+2},\dots,a_{\ell+\ell'-1}.
\end{equation*}
Here, $a_i$ ($1\leq i \leq \ell+\ell'-1$) are again some distinct
colours not contained in $\im c$. Any colouring $c'$ obtained in this
way from $c$ is said to be a \emph{continuation of $c$ to $H\cup Q\cup
  Q'$}. Note that $c'$ uses $\ell+\ell'-1$ colours not contained in
$\im c$, which is half the number of vertices in $V(Q\cup Q')\sm V(H)$.

\begin{figure}
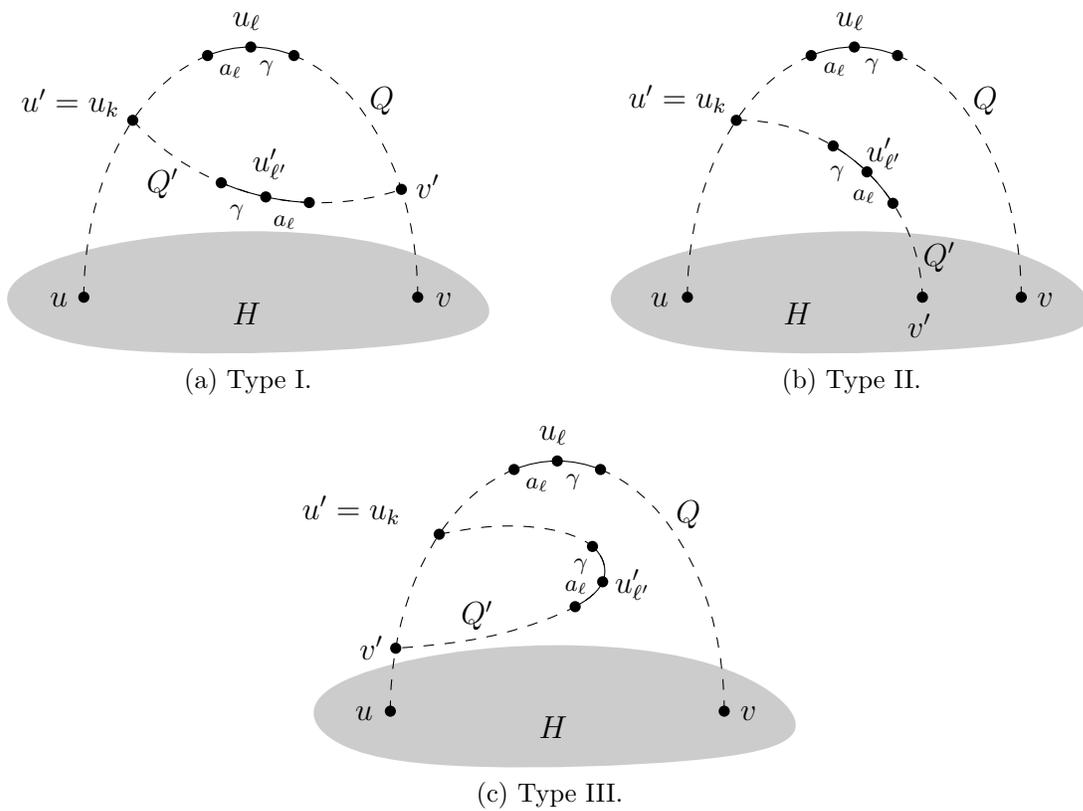

  \centering
  \sfig1{Type I.}\hf\sfig2{Type II.}\\
  \sfig3{Type III.}
  \caption{The possible types of the 2-extension $H\cup Q\cup
    Q'$. Dashed and solid lines represent paths and edges, respectively.}
  \label{fig:types}
\end{figure}

Based on the position of the endvertex $v'$ of $Q'$ relative to $Q$,
we distinguish three possible types of the 2-extension $H\cup Q\cup
Q'$. As shown in Figure~\ref{fig:types}, we may have $v'\in V(u'Qv)$
(type I), $v'\in V(H)\sm V(Q)$ (type II) or $v'\in V(uQu')$ (type
III). Note that types I and III include the possibility that $v'$
coincides with $v$ or $u$, respectively. Observe also that if $u' =
u_\ell$, then the 2-extension is of type I or II.

\begin{lemma}\label{l:two-even}
  Let $c$ be a safely rainbow-connecting colouring of $H$. If
  $H'=H\cup Q\cup Q'$ is an even 2-extension of $H$, then any
  continuation of $c$ to $H'$ is safely rainbow-connecting.
\end{lemma}
\begin{proof}
  Let $c'$ be a continuation of $c$ to $H'$. We use the same notation
  as in the definition of $c'$. We define $A$ as the set of colours
  used by $c'$ on $H\cup Q$.

  First, we show that $c'$ is rainbow-connecting. Let $x,y\in
  V(H')$. We are looking for a rainbow $xy$-path $R_{xy}$. If $x,y\in
  V(H\cup Q)$, then the argument is similar to that used in the proof
  of Lemma~\ref{l:odd}. For $x,y\in V(H)$, $R_{xy}$ is a rainbow
  $xy$-path in $(H,c)$. If $x\in V(Q)$ and $y\in V(H)$, then consider
  a $\gamma$-free subpath of $Q$ from $x$ to a vertex $w\in\Setx{u,v}$
  and define $R_{xy} = xQwR_{wy}y$ (where $R_{wy}$ has been defined
  before as $w\in V(H)$). Finally, if both $x$ and $y$ are vertices of
  $Q$, then we may assume that $x\in V(uQv)$; the path $R_{xy}$ is
  defined as $xQy$ if this path is rainbow, and $xQuR_{uv}vQy$
  otherwise.

  For later use, note that in all the cases considered up to now,
  $R_{xy}$ is either $a_\ell$-free or disjoint from $H$, with one
  exception, namely if $x = u_\ell$ and $y\in V(H)$.

  It remains to discuss the case that $x$ or $y$ is in $Q'$. By
  symmetry, we may assume that $x\in V(Q')$.

  \begin{xcase}{$y\in V(Q')$.}%
    If the path $xQ'y$ is not rainbow, then without loss of
    generality, we can write $x=u'_i$, $y=u'_j$, where $j\geq i +
    \ell' + 2$. In particular, $i \leq \ell'-2$ and $j \geq \ell'+2$,
    so $c'$ uses colours $a_{\ell+1},\dots,a_{\ell+i}$ on the path
    $xQ'u'$ and a subset of the colours
    $a_{\ell+i+2},\dots,a_{\ell+\ell'-1}$ on $v'Q'y$. It follows that
    the path $xQ'u'R_{u'v'}v'Q'y'$ is rainbow.
  \end{xcase}

  \begin{xcase}{$y\in V(H\cup Q)$ and $x \neq u'_{\ell'}$.}%
    If $x \neq u'_{\ell'}$, then there is an $A$-free path $S$ from
    $x$ to a vertex $w$ of $H\cup Q$ (just take a suitable subpath of
    $Q'$), and the path $xSwR_{wy}y$ is rainbow.
  \end{xcase}

  We are left with the following last case:
  \begin{xcase}{$y\in V(H\cup Q)$ and $x = u'_{\ell'}$.}%
    The path $xQ'u'Qu$ is rainbow as the colours from $A$ used on it
    are, in the order from $x$ to $u$,
    \begin{equation*}
      \gamma,a_{\ell'+\ell-1},a_{\ell'+\ell-2},\dots,a_{\ell+1},a_k,
      a_{k-1},\dots,a_1.
    \end{equation*}
    If $y$ is contained in this path, then the appropriate subpath is
    a rainbow $xy$-path. Similarly, if $y\in V(u_kQu_\ell)$, then there
    is a rainbow $xy$-path as $xQ'u'Qu_\ell$ is rainbow.

    The path $xQ'v'$ is rainbow and $c'$ uses no colours from $A$ on
    it except $a_\ell$. It follows that if $y\in V(H)$, then we can
    append either $v'QvR_{vy}y$ (for type I), $v'R_{v'y}y$ (for type
    II) or $v'QuR_{uy}y$ (for type III), and get a rainbow $xy$-path.

    It remains to consider the case that $y\in V(u_{\ell+1} Q
    v)$. Observe that the following subgraphs are rainbow and $A$-free
    under $c'$: $xQ'v' \cup u_{\ell+1}Qv$ (types I and II) and
    $xQ'v'Qu \cup u_{\ell+1}Qv$ (type III). Adding the rainbow path
    $R_{v'v}$ (type II) or $R_{uv}$ (type III) if necessary, we obtain
    a rainbow $xy$-path $R_{xy}$ in each of the cases.
  \end{xcase}
  
  It remains to check that $c'$ is safe. Observe first that whenever
  the above-defined path $R_{xy}$ is either $a_\ell$-free or
  edge-disjoint from $H$, then it is non-blocking (in the latter case,
  this is because $\size{\im c}\geq 2$, and only one colour from $\im
  c$ is used on $Q\cup Q'$ by $c'$).

  By directly inspecting the above construction, we can readily check
  that the following cases are (up to symmetry) the only ones where
  $R_{xy}$ is neither $a_\ell$-free nor edge-disjoint from $H$:
  \begin{enumerate}[\quad(a)]
  \item $x=u_\ell$ and $y\in V(H)$,
  \item $x=u'_{\ell'}$ and $y\in V(H)$,
  \item $x=u'_{\ell'}$, $y\in V(u_{\ell+1}Qv)$ and the 2-extension $H\cup
    Q\cup Q'$ is of type II or III.
  \end{enumerate}

  We will now show that only at most one vertex is blocked for
  $u_\ell$. If $y$ is such a vertex, then we are in case (a) above. By
  the construction, $R_{xy} = xQuR_{uy}y$, so $R_{uy}$ must be
  blocking in $(H,c)$. By the assumption that $c$ is safe, there is at
  most one such vertex $y$ as claimed.

  Next, we consider the vertex $x=u'_{\ell'}$. Suppose that some
  vertex $y\in V(u_{\ell+1}Qv)$ is blocked for $x$ (case (c)). Since
  the 2-extension must be of type II or III, we have $y=u_{\ell+1}$,
  for otherwise the colour $a_1$ is not used on $R_{xy}$. Furthermore,
  no vertex $y'\in V(H)$ is blocked for $y$ as $a_1$ is not used on
  $R_{xy'}$. Thus, $u_{\ell+1}$ is the only blocked vertex for
  $u'_{\ell'}$.

  For this choice of $x$, it remains to consider the case that no
  vertex of $u_{\ell+1}Qv$ is blocked for $x$. Suppose that $y\in V(H)$
  is blocked for $x$ (case (b)). By the construction, $R_{xy}$
  contains as a subpath the path $R_{vy}$ (for type I), $R_{v'y}$
  (type II) or $R_{uy}$ (type III). In addition, all the colours from
  $\im c$ are used on this subpath. It follows that in $(H,c)$, $y$ is
  blocked for $v$, $v'$ or $u$ depending on the type, so $y$ is
  uniquely determined since $c$ is safe.

  We have shown that for $x\in\Setx{u_\ell,u'_{\ell'}}$, there is at
  most one $y$ which is blocked for $x$. Considering the cases
  (a)--(c) above, $c'$ will be proved safe if we show that no $y\in
  V(H)$ is blocked for both $u_\ell$ and $u'_{\ell'}$. Suppose that
  $y\in V(H)$ is blocked for both of these vertices. Since $R_{u_\ell
    y} = u_\ell QuR_{uy}y$, all the colours from $\im c$ must be used
  on $R_{uy}$, so $u$ is blocked for $y$ in $(H,c)$. By similarly
  considering $R_{u'_{\ell'}y}$, we find that for type I or II, the
  vertex $v$ or $v'$, respectively, would also be blocked for $y$ in
  $(H,c)$, which is impossible as $c$ is safe.

  Hence, the 2-extension must be of type III. In this case, as
  observed in the definition of types, $u' \neq u_\ell$ and hence $v'
  = u_i$ with $i \leq \ell-2$. Consequently, the colour $a_{\ell-1}$
  is not used by $c'$ on $R_{u'_{\ell'}y}$, contradicting the
  assumption that $y$ is blocked for $u'_{\ell'}$. The proof is
  complete.
\end{proof}


\section{Proof of Theorem~\ref{t:main}}
\label{sec:proof-theorem}

In this section, we prove Theorem~\ref{t:main}. Let $G$ be a
2-connected graph with $n$ vertices. 

If $G$ is an odd cycle, then the colouring of its edges by
\begin{equation*}
  a_1,a_2,\dots,\lfloor n/2\rfloor,\lceil n/2\rceil,1,2,\dots,\lfloor
  n/2\rfloor
\end{equation*}
is a rainbow colouring with $\lceil n/2\rceil$ colours. Thus, we may
assume that $G$ is not an odd cycle.

We claim that $G$ contains an even cycle. If not, let $Z$ be an odd
cycle in $G$ (recall that $G$ is 2-connected) and let $v$ be a vertex
not contained in $Z$. Taking two internally disjoint paths from $v$ to
distinct vertices $z_1,z_2$ on $Z$ and concatenating them with a
$z_1z_2$-subpath of $Z$ with the appropriate parity, we obtain an even
cycle, a contradiction.

Thus, let $H_0$ be an even cycle in $G$, say of length $2k$. We
construct a subgraph $\Hstar$ of $G$ by means of a sequence $H_0, H_1,
\dots$ of subgraphs of $G$. To construct $H_{i+1}$ ($i\geq 0$), we
proceed as follows:
\begin{itemize}
\item if there is an odd $H_i$-path $P_i$, we set $H_{i+1} = H_i \cup
  P_i$ (making an arbitrary choice if there are more such paths),
\item otherwise, if there is a 2-extension $H_i\cup Q_i\cup Q'_i$ of
  $H_i$ with both $Q_i$ and $Q'_i$ even, we set $H_{i+1} = H_i\cup
  Q_i\cup Q'_i$,
\item if there is neither an odd $H_i$-path nor an even 2-extension of
  $H_i$, we finish and set $\Hstar = H_i$.
\end{itemize}

In the rest of this section, the symbol $\Hstar$ will denote the
subgraph of $G$ just constructed. Observe that in the above sequence,
each subgraph $H_i$ has an even number of vertices. Thus,
$\size\Hstar$ is even. The following proposition describes the weak
$\Hstar$-bridges.

\begin{proposition}\label{p:bridges}
  Let $B$ be a weak $\Hstar$-bridge. Then the following holds:
  \begin{enumerate}[\quad(i)]
  \item $\Hstar\cup B$ is an extension of $\Hstar$ by a path sequence
    $(P_1,P_2,\dots,P_\ell)$, where $P_1$ is even and the other paths
    are odd,
  \item $\size B$ is odd.
  \end{enumerate}
\end{proposition}
\begin{proof}
  (i) Let $M$ be an inclusionwise maximal extension of $\Hstar$
  contained in $\Hstar\cup B$. Choose a path sequence
  $(Q_1,\dots,Q_s)$ for $M$. Note that by the construction of $\Hstar$ and
  Lemma~\ref{l:paths}, $Q_1$ must be even and all the other paths
  $Q_i$ must be odd.

  We claim that $M=\Hstar\cup B$. Suppose that this is not the case and
  choose a vertex $w\in V(B)\sm V(M)$. By the 2-connectedness of $G$,
  there are internally disjoint paths $R_1,R_2$ from $w$ to distinct
  vertices of $M$. The concatenation of $R_1$ and $R_2$ is an $M$-path
  of length at least 2 which can be added to $M$ and provides a
  contradiction with the maximality of $M$. This proves part (i).

  Part (ii) is a direct consequence of (i).
\end{proof}

On each $H_i$, we define a safely rainbow-connecting colouring $c_i$
by $\size{H_i}/2$ colours. We begin with a colouring $c_0$ of the even
cycle $H_0$ with values
\begin{equation*}
  1,2,\dots,k,1,2,\dots,k.  
\end{equation*}
It is easy to check that this colouring is safely
rainbow-connecting. Given $c_i$, the colouring $c_{i+1}$ is
constructed as a continuation to $H_i\cup P_i$ or $H_i\cup Q_i\cup
Q'_i$, respectively. By Lemmas~\ref{l:odd} and \ref{l:two-even}, we
eventually obtain a safely rainbow-connecting colouring $\cstar$ of
$\Hstar$. By the construction, $\cstar$ uses $\size \Hstar/2$ colours.

At this point, we fix two more `special' colours in addition to
$\gamma$, namely $\alpha$ and $\beta$. The choice is such that neither
$\alpha$ nor $\beta$ is contained in $\im\cstar$.

Let $u,v\in V(\Hstar)$ and let $H$ be a subgraph of $G$ such that
$\Hstar\subseteq H$. A colouring $b$ of $H$ is \emph{grounded in
  $(u,v)$} if $b$ coincides with $\cstar$ on $\Hstar$, and for each
vertex $x\in V(H)\sm V(\Hstar)$, both of the following conditions hold:
\begin{enumerate}[\quad (\textrm{A}1)]
\item $(H-E(\Hstar),b)$ contains either a rainbow
  $\alpha\beta\gamma$-free path from $x$ to $\Hstar$, or both a
  rainbow $\beta\gamma$-free $xu$-path and a rainbow
  $\alpha\gamma$-free $xv$-path,
\item for every vertex $y\in V(H)\sm V(\Hstar)$, there is a rainbow
  $xy$-path in $H$ which is either edge-disjoint from $\Hstar$, or
  $\beta$-free.
\end{enumerate}
Note that the definition of a grounded colouring is always related to
the same subgraph $\Hstar$ of $G$ defined above. If the pair $(u,v)$
is not essential, then we just say that $b$ is \emph{grounded}.

The following lemma shows that a continuation of a grounded colouring
to an odd 1-extension is grounded.
\begin{lemma}\label{l:odd-grounded}
  Let $u,v\in V(\Hstar)$ and let $c$ be a colouring of $H$,
  $\Hstar\subseteq H\subseteq G$, which is grounded in $(u,v)$. If $P$ is
  an $H$-path of odd length, then any continuation of $c$ to $H\cup
  P$ is grounded in $(u,v)$.
\end{lemma}
\begin{proof}
  Let $c'$ be a continuation of $c$ to $H\cup P$ and let $x$ be a
  vertex in $V(H\cup P)\sm V(\Hstar)$. It suffices to verify properties
  (A1) and (A2) for $x\in V(P)$ since for any other $x$ they follow
  from the assumption.

  We begin with (A1). Let $P'$ be the shorter of the two subpaths of
  $P$ with one endvertex $x$ and the other endvertex in $H$. Let the
  latter endvertex be denoted by $w$. Observe that $P'$ is rainbow and
  $\alpha\beta\gamma$-free in $(H\cup P,c')$. By the assumption on
  $c$, $(H,c)$ contains either a rainbow $\alpha\beta\gamma$-free
  path $R$ from $w$ to a vertex $z\in V(H)$, or a $\beta\gamma$-free
  $wu$-path $R_1$ and an $\alpha\gamma$-free $wv$-path $R_2$, both
  rainbow. In the former case, the path $xP'wRz$ is rainbow and
  $\alpha\beta\gamma$-free, because $P'$ is $\alpha\beta\gamma$-free
  and no colour from $\im c$ is used on $P'$. In the latter case, we
  similarly obtain paths with the desired properties by prepending
  $P'$ to $R_1$ and $R_2$, respectively.

  To verify (A2), let $y$ be a vertex in $V(H\cup P)\sm V(\Hstar)$. If
  $y\notin V(P)$, then we can utilise the path $P'$ as above and
  concatenate it with a rainbow $wy$-path in $(H,c)$ satisfying
  (A2); the resulting $xy$-path satisfies (A2) as well, because $P'$
  is $\beta\gamma$-free and edge-disjoint from $\Hstar$.

  It remains to discuss the case that $y\in V(P)$. If $xPy$ is
  rainbow, then we are done as it is $\beta$-free. Otherwise, let the
  endvertices of $P$ be denoted by $x'$ and $y'$ in such a way that
  $x\in V(x'Py)$. Observe that $xPx'\cup yPy'$ is rainbow and
  $\alpha\beta\gamma$-free.

  If the vertices $x',y'$ are not in $\Hstar$, then by the assumption
  that $c$ is grounded, a rainbow $x'y'$-path $S$ in $(H,c)$ is either
  edge-disjoint from $\Hstar$, or $\beta$-free. It follows that the
  $xy$-path $xPx'Sy'Py$ is edge-disjoint from $\Hstar$ or $\beta$-free
  as well. Since $\gamma$ may be used on $S$, it is important that
  $xPx' \cup yPy'$ is $\gamma$-free. This makes the $xy$-path rainbow.

  We may therefore assume that $y'\in V(\Hstar)$. If $x'\in
  V(\Hstar)$, then since $\cstar$ is rainbow-connecting, there is a
  rainbow $x'y'$-path $S_1$ in $(\Hstar,\cstar)$ which is $\beta$-free
  as $\beta\notin\im\cstar$. Consequently, the $xy$-path $xPx'Sy'Py$ is
  rainbow and $\beta$-free.

  Lastly, if $x'\in V(H)\sm V(\Hstar)$ and $y'\in V(\Hstar)$, then by the
  assumption that $c$ is grounded, $H-E(\Hstar)$ contains a rainbow
  $\beta\gamma$-free path $S_2$ from $x'$ to a vertex $z\in
  V(\Hstar)$. Furthermore, there is a rainbow (necessarily
  $\beta$-free) $zy'$ path $S_3$ in $(\Hstar,\cstar)$. The path
  $xPx'S_2zS_3y'Py$ is then $\beta$-free and rainbow.
\end{proof}

Let $B^1,\dots,B^r$ be all the weak $\Hstar$-bridges in $G$. Fix $j$,
$1\leq j \leq r$. We use Proposition~\ref{p:bridges} to choose a path
sequence $(P^j_1,\dots,P^j_{\ell^j})$ for $\Hstar\cup B^j$. The only
even path in this sequence, $P^j_1$, will be called the \emph{base} of
$B^j$. We let the length of $P^j_i$ be denoted by $2k^j_i$ or
$2k^j_i+1$, according to whether it is even or odd. Furthermore, the
endvertices of the base of $B^j$ will be denoted by $u^j$ and $v^j$.

We now extend $\cstar$ to a suitable colouring $c^j$ of $\Hstar\cup
B^j$. We first colour the base of $B^j$, in the direction from $u^j$
to $v^j$, by
\begin{equation*}
  a_1,a_2,\dots,a_{k^j_1-1},\alpha,\beta,a_1,a_2,\dots,a_{k^j_1-1},
\end{equation*}
where the $a_i$ ($1\leq i \leq k^j_1-1$) are distinct colours not
contained in $\im\cstar\cup\Setx{\alpha,\beta}$ nor used for the
colouring of any other weak $\Hstar$-bridge. The colours $\alpha$ and
$\beta$ are used for all the bases. Observe that the colouring is
grounded in $(u^j,v^j)$.

We extend the colouring to all of $\Hstar\cup B^j$ by successively
taking continuations to odd 1-extensions in the above path sequence
for $\Hstar\cup B^j$ (recall that all the paths $P^j_i$ with $i\geq 2$
are odd). By a repeated use of Lemmas~\ref{l:odd} and
\ref{l:odd-grounded}, we obtain a safely rainbow-connecting colouring
$c^j$ which is grounded in $(u^j,v^j)$. Let $x^j$ be the number of
colours used by $c^j$ on $B^j-E(\Hstar)$. We compare $x^j$ to the
number of vertices in $V(B^j)\sm V(\Hstar)$. For $2\leq i\leq \ell^j$,
the path $P^j_i$ has length $2k^j_i+1$, and therefore $2k^j_i$
internal vertices. The number of vertices in $V(B^j)\sm V(\Hstar)$ is thus
$2k^j_1 - 1 + 2k^j_2 + \dots + 2k^j_{\ell^j}$. On the other hand,
\begin{equation}\label{eq:xj}
  x^j = k^j_1 + 1 +
  \sum_{i=2}^{\ell^j} k^j_i =  \frac{\size{V(B^j)\sm V(\Hstar)}}2 + \frac32.
\end{equation}

Since each $c^j$ extends $\cstar$ and the weak $\Hstar$-bridges are
pairwise edge-disjoint, we can combine all the colourings $c^j$,
$1\leq j\leq r$, to a colouring $\tilde c$ of $G$. We assert that
$(G,\tilde c)$ is rainbow-connected. To check this, it is enough to
show that any two vertices contained in different $\Hstar$-bridges are
joined by a rainbow path. Let us say that $x\in V(B^1)\sm V(\Hstar)$ and
$y\in V(B^2)\sm V(\Hstar)$. By condition (A1) in the definition of a
grounded colouring, $(\Hstar\cup B^1,c^1)$ contains an
$\alpha\gamma$-free rainbow path from $x$ to a vertex $w_1$ in
$\Hstar$, and $(\Hstar\cup B^2,c^2)$ contains a $\beta\gamma$-free
rainbow path from $y$ to a vertex $w_2$ in $\Hstar$. Since no colour
from $\im\cstar$ is used on these paths by $\tilde c$, and each of
$\alpha$ and $\beta$ is used on at most one of them, we can
concatenate the paths with a rainbow $w_1w_2$-path in
$(\Hstar,\cstar)$ and obtain a rainbow $xy$-path in $(G,\tilde
c)$. This completes the proof that $(G,\tilde c)$ is rainbow.

The number of colours used by $\tilde c$ can be obtained by
using~\eqref{eq:xj} and making a correction to account for the fact
that the same colours $\alpha,\beta$ are used in all $B^j$:
\begin{align}\label{eq:ncolours}
  \size{\im{\tilde c}} &= \size{\im\cstar} +
  \Bigl(\sum_{j=1}^r x^j\Bigr) - 2(r-1)\notag\\
  &= \frac{\size\Hstar}2 + \Bigl(\sum_{j=1}^r
  \frac{\size{V(B^j)\sm V(\Hstar)}}2 \Bigr) - \frac r2 + 2 = \frac{\size
    G}2 - \frac r2 + 2.
\end{align}
Note that for $r \geq 3$, we have $\size{\im{\tilde c}} \leq \lceil
\size G/2\rceil$, which implies the statement in
Theorem~\ref{t:main}. Similarly, if $r=0$, then $\Hstar=G$ and we are
done as well, because $(G,\cstar)$ is then rainbow-connected and
$\cstar$ uses $\size G/2$ colours. Consequently, we may assume that
\begin{equation*}
  1\leq r\leq 2.
\end{equation*}
We will perform a simple recolouring to reduce the number of colours
by one and obtain a colouring satisfying the bound in
Theorem~\ref{t:main}.
 
The case $r=1$ is simple. Since $\tilde c$ is grounded, every vertex
in $V(G)\sm V(\Hstar)$ is joined to $\Hstar$ by a $\beta\gamma$-free
path. Thus, if we recolour the unique edge coloured by $\beta$ to a
colour $\gamma'\in\im\cstar\sm \Setx{\gamma}$, $G$ is still
rainbow-connected with respect to the resulting colouring. The number
of colours used after this reduction is exactly $\lceil\size G/2\rceil$.

It thus remains to consider the case $r=2$. This case is the reason
why we need to restrict ourselves to safe colourings. For $j=1,2$, let
$e^j_\alpha$ and $e^j_\beta$ denote the edge of $P^j_1$ coloured by
$\alpha$ and $\beta$, respectively.

Since $\cstar$ is safe, there is a rainbow path $P$ from $u^1$ to
either $u^2$ or $v^2$, and a colour $\delta\in\im\cstar$ which is not
used on $P$. (We allow $\delta = \gamma$.) As the reversal of the
colouring on $P^2_1$ affects neither the rainbow-connectedness of $G$
nor the groundedness of the colouring, we may actually assume that the
endvertices of $P$ are $u^1$ and $v^2$.

We recolour the edges $e^1_\beta$ and $e^2_\beta$ to $\delta$, and we
argue that with the resulting colouring $\bar c$, $G$ is still
rainbow-connected. Let $x,y$ be vertices of $G$. Since no rainbow
paths inside $\Hstar$ are affected by the recolouring, we may assume
that $x\in V(B^1) \sm V(\Hstar)$.

We distinguish several cases based on the location of $y$. 

\setcounter{xcasehdr}0
\begin{xcase}{$y\in V(B^1)\sm V(\Hstar)$.}%
  Since $c^1$ is grounded, condition (A2) says that in $(\Hstar\cup
  B^1,c^1)$ there is a rainbow $xy$-path which is either $\beta$-free,
  or edge-disjoint from $\Hstar$. In either case, the colouring $\bar c$
  only uses the colour $\delta$ at most once on this path, which is
  therefore rainbow.
\end{xcase}

\begin{xcase}{$y\in V(\Hstar)$.}%
  Since $c^1$ is grounded, $(B^1,c^1)$ contains a $\beta\gamma$-free
  rainbow path $Q_1$ from $x$ to a vertex $w\in V(\Hstar)$. If we let
  $Q_2$ be a rainbow $wy$-path in $(\Hstar,\cstar)$, then $xQ_1wQ_2y$ is a
  rainbow $xy$-path.
\end{xcase}

\begin{xcase}{$y\in V(B^2)\sm V(\Hstar)$.}%
  As $c^1$ is grounded, there is a rainbow path $R_1$ in $(B^1,c^1)$
  from $y$ to a vertex $w_1\in V(\Hstar)$, where $R_1$ is either
  $\alpha\beta\gamma$-free, or it is $\beta\gamma$-free and $w_1=u^1$.

  As a first subcase, assume that $R_1$ is $\alpha\beta\gamma$-free,
  and choose a rainbow $\beta\gamma$-free path $R_2$ in $(B^2,c^2)$
  from $y$ to a vertex $w_2\in V(\Hstar)$ using condition (A1) in the
  definition of a grounded colouring. Thus, $R_1\cup R_2$ is
  $\beta\gamma$-free with respect to $\tilde c$. Since in $(G,\tilde
  c)$, the sets of colours used to colour $B^1$ and $B^2$ are disjoint
  except for $\alpha,\beta,\gamma$, $R_1\cup R_2$ is rainbow with
  respect to $\tilde c$. 

  Let $R_0$ be a rainbow $w_1w_2$-path in $(\Hstar,\cstar)$, and define
  $R=xR_1w_1Rw_2R_2y$. Since $\gamma$ is not used by $\tilde c$ on
  $R_1\cup R_2$, $R$ is rainbow in $(G,\tilde c)$. Furthermore, it
  remains rainbow after the recolouring of $e^1_\beta$ and $e^2_\beta$
  to $\delta$, since $R_1\cup R_2$ is $\beta$-free. Thus, $R$ is
  rainbow in $(G,\bar c)$.

  A symmetric situation occurs when $(B^2,c^2)$ contains an
  $\alpha\beta\gamma$-free rainbow path from $y$ to $\Hstar$. We may
  therefore assume that $R_1$ is a $\beta\gamma$-free $xu^1$-path, and
  choose an $\alpha\gamma$-free rainbow $yv^2$-path $R'_2$ in
  $(\Hstar\cup B^2,c^2)$. Let $P$ be the $\delta$-free rainbow
  $u^1v^2$-path in $(\Hstar,\cstar)$ defined above, and $R' =
  xR_1u^1Pv^2R'_2y$.

  We claim that no colour is repeated on $R'$ in $(G,\bar c)$:
  $\alpha$ may only used on $R_1$ as $R'_2$ is $\alpha$-free and
  $\alpha\notin\im\cstar$, $\gamma$ is not used on $R_1\cup R'_2$, and
  $\delta$ may only be used on $R'_2$, since $P$ is $\delta$-free and
  so is $R_1$ (being chosen to be $\beta$-free with respect to
  $c^1$). This concludes the discussion of this case.
\end{xcase}

We have just shown that $(G,\bar c)$ is rainbow-connected. As for
$\size{\im{\bar c}}$, by eliminating the color $\beta$ we decreased
the value in \eqref{eq:ncolours} by one, and we get
\begin{equation*}
  \size{\im{\bar c}} = \frac{\size G}2.
\end{equation*}
The proof of Theorem~\ref{t:main} is now complete.


\section{Higher connectivity}
\label{sec:connectivity}

In view of Theorem~\ref{t:main}, it is natural to ask whether one can
further improve the bound for graphs of higher connectivity. The
question on the relation between connectivity and the rainbow
connection number was previously asked by H.~J.~Broersma
(cf.~\cite[Problem 2.22]{LS-survey}). The best bound to date follows
from a result of Chandran et al.~\cite{CDRV}, improving a bound of
Schiermeyer~\cite{S2}:
\begin{theorem}[\cite{CDRV}]\label{t:degree}
  A connected graph $G$ with $n$ vertices and minimum degree
  $\delta(G)$ has
  \begin{equation*}
    \rconn G \leq \frac{3n}{\delta(G)+1} + 3.
  \end{equation*}
\end{theorem}
Since the minimum degree of a graph is greater than or equal to its
connectivity, the theorem implies a bound for $k$-connected graphs.

\begin{figure}
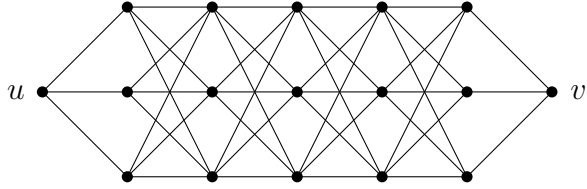

  \centering
  \fig4
  \caption{A $k$-connected graph $G$ with $k(\ell-1)+2$ vertices and
    $\rconn G \geq \ell$. (Shown for $k=3$ and $\ell=6$.)}
  \label{fig:lower}
\end{figure}

We show that for every $k$, there are $k$-connected graphs $G$ with
$n$ vertices and
\begin{equation*}
  \rconn G \geq \frac{n-2}k + 1.
\end{equation*}
For fixed $k,\ell$, let $P$ be a path of length $\ell$ with
endvertices $u_0$ and $v_0$, and let $I$ be a graph consisting of $k$
independent vertices. Let $G_0$ be the lexicographic product of $P$
and $I$. Thus, $G_0$ has vertex set $V(P)\times V(I)$ and vertices
$(x,y)$ and $(x',y')$ are joined by an edge whenever $x$ and $x'$ are
adjacent in $P$. Let $G$ be the graph obtained from $G_0$ by
identifying all vertices of the form $(u_0,y)$ ($y\in V(I)$) into one
vertex $u$, and all vertices of the form $(v_0,y)$ ($y\in V(I)$) into
another vertex $v$. (See Figure~\ref{fig:lower}.)

It is easy to see that $G$ is $k$-connected and has $n := k(\ell-1)+2$
vertices. Since shortest $uv$-paths have length $\ell$, we have
\begin{equation*}
  \rconn G \geq \ell = \frac{n-2}k + 1
\end{equation*}
as claimed above. We conclude with a question:
\begin{problem}
  Is there a constant $C=C(k)$ such that every $k$-connected graph $G$
  with $n$ vertices satisfies
  \begin{equation*}
    \rconn G \leq \frac nk + C?
  \end{equation*}
\end{problem}


\end{document}